\documentclass[12pt]{article}

\usepackage[utf8]{inputenc}
\usepackage[margin=3cm]{geometry}
\usepackage{amsfonts, amsmath, amssymb, amsthm}
\usepackage{indentfirst}

\newtheorem{Dfn}{Definition}
\newtheorem{Thm}{Theorem}
\newtheorem{Prp}{Proposition}

\newenvironment{Prf}{\par\noindent{\it Proof:}\rm}{\newline\rightline{\textbf{QED}}\par}
\newtheorem{Cor}{Corrolary}

\newcommand{\R}{\mathbb{R}}

\newcommand{\con}{\mathfrak{c}}

\newcommand{\Cl}{\rm{Cl}}
\newcommand{\Int}{\rm{Int}}
\newcommand{\Pow}{\mathcal{P}}

\newcommand{\sat}{\rm{sat}}

\newcommand{\dom}{\mathrm{dom}}

\title{Connections between Kuratowski partitions of Baire spaces, measurable cardinals and precipitous ideals}

\author{Sławomir Kusiński}

\date{}

\begin{document}
    \maketitle
    
    \begin{abstract}
        In this paper we present a few properties of \(K\)-partitions, which are partitions of Baire spaces such that all subfamilies of such a partition sum to a set with the Baire property. Among the result proven we have general existence result that state that the existence of any \(K\)-partition implies the existence of \(K\)-partition of a metrizable space as well as existence of \(K\)-partition of a compact space implies the existence of \(K\)-partition of a completely metrizable space. We also prove some connections between existence of \(K\)-partitions and existence of precipitous ideals as well as measurable cardinals. There are also outlined possible connection with real-measurable cardinals, extensions of Lebesgue measure on the closed interval and density topologies. 
    \end{abstract}
    
    \section{Introduction}
    
    In this paper we take a look at connections between Kuratowski partitions of Baire spaces with other foundational concepts such as measurable cardinals and precipitous ideals. Kuratowski partition or shortly \(K\)-partitions are a type of very well behaved partitions on Baire spaces. The caveat is that do not need to exist under ZFC, even more so their existence is not equiconsistent with ZFC. Nevertheless their deep connections with other mathematical subjects makes studying them worthwhile.
    
    Their genesis traces back to the paper \cite{Kur} by Kuratowski where he posed a question about functions with the Baire property from a completely metrizable space to a metrizable space, exactly do and when such functions have to be continuous apart from a meager set or not. From the beginning it was known that it was the case when the range of a function was separable. The real question was when the assumption of separability could be dropped. As was later shown in \cite{Fr} such a question could be rephrased using \(K\)-partitions. One of the earliest approaches to the subject was by Kunugi\cite{Kunugi}, a fact that has been lately reviewed in \cite{GrzeLab}.
    
    \section{Definitions and known facts}
    
    Throughout the whole paper we assume that we are dealing with Hausdorff topological spaces exclusively and the space denoted by \(X\) is assumed to be Baire. By a Baire space we mean a topological space in which any countable intersection of open dense sets is not empty, ie a space that is non-meager in itself. Some authors would call that a localy Baire space. 
    
    \begin{Dfn}
        Let \( X \) be a topological space and \( A \subseteq X \). \( A \) is said to have Baire property if it can be represented as \( U \triangle F \), where \( U \) is open, \( F \) is meager and \( \triangle \) denotes the symmetric difference of sets.
    \end{Dfn}

    \begin{Dfn}
        Let \(X\) be a topological space. A partition \( \mathcal{F} \) of \( X \) into meager sets is called a \(K\)-partition if for any \( \mathcal{F}' \subseteq \mathcal{F} \) the set \( \bigcup \mathcal{F}' \) has Baire property.
    \end{Dfn}

    \begin{Dfn}
        Let \(Y\) be a boolean algebra and \( I \) be an ideal on \( Y \). The saturation of \( I \) - denoted by \( \sat(I) \) - is the smallest cardinal such that all antichains in \( Y / I \) are of cardinality less than \( \sat(I) \).
    \end{Dfn}

    \begin{Dfn}
        Let \( I \) be an ideal on a set \( Y \). The set \( W \subseteq \Pow(Y) \) is called an \( I \)-partition of \(Y\) if \( \bigcup W = Y \) and for any \( A_1, A_2 \in W \) if \( A_1 \not= A_2 \) then \( A_1 \cap A_2 \in I \). The ideal \( I \) is called precipitous if it is \( |Y| \)-complete and for any sequence \( W_0, W_1, \ldots \) of \( I \)-partitions such that \( W_{n+1} \) is a refinement of \( W_n \) there exist \( X_i \in W_i \) such that \( \bigcap\limits_{i \in \omega} X_i \not= \emptyset \).
    \end{Dfn}

    For any ideal \( I \) on a cardinal \( \kappa \) we will denote \( I^+ = P(\kappa) \setminus I \). Now let
    \[
        X(I) = \{ x \in (I^+)^\omega \colon \forall_{n \in \omega} \bigcap\limits_{k=0}^n x(k) \in I^+, \, \bigcap\limits_{k < \omega} x(k) \not= \emptyset \}
    \]
    be a subspace of the metric space \( (I^+)^\omega \) where \( I^+ \) is a discrete space. It has been shown in \cite{KunFr} that \( I \) is precipitous iff \( X(I) \) is a Baire space and if \( I \) is precipitous then sets
    \[
        F_\alpha = \{ x \in X(I) \colon \alpha = \min \bigcap\limits_{k < \omega} x(k) \} \text{ for } \alpha \in \kappa
    \]
    define a \( K \)-partition of \( X(I) \).
    
    Let \( \tau^+ \) be a discrete space of all non empty open subsets of \( X \). Similarly we can define
    \[
        X(\tau) = \{ x \in (\tau^+)^\omega \colon \bigcap\limits_{k < \omega} x(k) \not= \emptyset \}
    \]
    and
    \[
        X^*(\tau) = \{ x \in (\tau^+)^\omega \colon \bigcap\limits_{k < \omega} x(k) \not= \emptyset, \forall_{n \in \omega} \Cl(x(n+1)) \subseteq x(n) \}.
    \]
    As we will see those spaces will be vital in showing that if there exists a \(K\)-partition of any Baire space then there also exists a \(K\)-partition of some metric space. This in fact will give us equiconsistency of existence of a measurable cardinal and the existence of a \(K\)-partition of any Baire space as it was shown in \cite{KunFr} that:
    \begin{Thm}
        The following theories are equiconsistent:
        \begin{itemize}
            \item ZFC \(+\) existence of a measurable cardinal,
            \item ZFC \(+\) existence of a \(K\)-partition of a Baire metric space,
            \item ZFC \(+\) existence of a \(K\)-partition of a complete metric space.
        \end{itemize}
    \end{Thm}
    
    If \( \mathcal{F} = \{ F_\alpha \colon \alpha \in \kappa \} \) is a \(K\)-partition of Baire space \( X \) then the set
    \[ I_\mathcal{F} = \{ A \in P(\kappa) \colon \bigcup\limits_{\alpha \in A} F_\alpha \mbox{ is meager } \} \]
    is an ideal on \( \kappa \). If \( U \subseteq X \) is open and non-meager then
    \[ \mathcal{F}|_U = \{ F_\alpha \cap U \colon \alpha \in \kappa \} \]
    is a \(K\)-partition of \( U \). What's more if \( U \subseteq V \) then \( I_{\mathcal{F}|_V} \subseteq I_{\mathcal{F}|_U} \) and thus \( X(I_{\mathcal{F}|_U}) \subseteq X(I_{\mathcal{F}|_V}) \).

    There is yet another equivalent approach to precipitous ideals.
    \begin{Dfn}
        Let \( I \) be an ideal on a set \( Y \). Let \( F \) be a family of functions on subsets of \( Y \) to ordinals. We will say that \( F \) is a functional if \( W_F = \{ \dom(\varphi) \colon \varphi \in F \} \) is an \(I\)-partition of \( Y \). Furthermore for two functionals \( F, G \) we define \( F < G \) if
        \begin{itemize}
            \item \( W_F \) refines \( W_G \)
            \item for \(f \in F, g \in G \) such that \( \dom(f) \subseteq \dom(g) \) we have \( f(\alpha) < g(\alpha) \) for \( \alpha \in \dom(f) \)
        \end{itemize}
    \end{Dfn}
    In \cite{Jech} a following characterization of precipitousness was shown.
    \begin{Thm}
        An ideal \( I \) is precipitous if and only if for no \( S \in I^+ \) there exists a sequence of functionals \( F_0 > F_1 > \ldots \)
    \end{Thm}

    In \cite{KurTop} one can find a following result attributed to Banach, called localization theorem.
    \begin{Thm}
        Let \( X \) be a topological space and let the sets \( U_i \) for \( i \in \kappa \) be open meager subsets of \( X \). Then \( \bigcup\limits_{i \kappa} U_i \) is also meager.
    \end{Thm}
    The localization theorem has one important consequence, which can be stated as follows.
    \begin{Cor}
        Let \( X \) be a Hausdorff and Baire space. There exists open subset \( U \subseteq X \) such that \( U \) has no non-empty meager open subset.
    \end{Cor}
    \begin{Prf}
        By the localization theorem the set \( W = \bigcup \{ V \subseteq X \colon V \mbox{ is open and meager } \} \) is open and meager. \( W \) cannot be dense, because \( X \) is a Baire space. Then \( U = \Int (X \setminus W) \) is as required.
    \end{Prf}
    A space with such properties may be called globally Baire. Those are exactly the spaces in which a countable intersection of open and dense subsets is dense.

    \section{Basic results about K-partitions}
    
    We will begin by proving that the existence of any \(K\)-partition whatsoever implies the existence of \(K\)-partition of some metrizable space. 
    
    \begin{Thm}
        Let \(X\) be a space with \( K \)-partition \( \mathcal{F} \) of minimal cardinality \( \kappa \) and let \(\tau\) be the topology of \(X\). Then the space \(X(\tau)\) has a \(K\)-partition.
    \end{Thm}
    \begin{Prf}
        First we will show that \(X(\tau)\) is a Baire space. Let \(G_i \subseteq X(\tau) \) be open and dense for \( i \in \omega \). Let
        \[
            G_i^\# = \bigcup \{ U_0 \cap \ldots \cap U_n \colon n \in \omega, \forall_{x \in X(\tau)} ( (x_0 = U_0 \wedge \ldots \wedge x_n = U_n) \Rightarrow x \in G_i ) \}.
        \] 
        Each \( G_i^\# \) is a sum of open sets an thus is open in \(X\). Let \( U \subseteq X \) be open and 
        \[
            G_U = \{ x \in X(\tau) \colon x_0 = U \}.
        \]
        \( G_U \) is open in \(X(\tau)\) and by density of \( G_i \) we get that \( G_U \cap G_i \) is non-empty and open. This means that there exist \( U_1, \ldots, U_n \) such that \( U \cap U_1 \cap \ldots \cap U_n \subseteq G_i^\# \) and thus \( U \cap G_i^\# \not= \emptyset \) which proves that \( G_i^\# \) is dense. Then as \(X\) is a Baire space we have \( y \in \bigcap\limits_{i \in \omega} G_i^\#  \) ie \( y \in U_{i,1} \cap \ldots \cap U_{i,n_i} \subseteq G_i^\# \). It follows that \( x = (U_{1,1}, \ldots, U_{1,n_1}, U_{2,1}, \ldots) \in \bigcap\limits_{i \in \omega} G_i \) and thus \( X(\tau) \) is a Baire space as required.
        
        Let
        \[
            \tilde{F}_\alpha = \{ x \in X(\tau) \colon \alpha = \min \{ \beta < \kappa \colon \bigcap\limits_{k < \omega} x(k) \cap F_\beta \not= \emptyset \} \} \text{ for } \alpha \in \kappa.
        \]
        We will show that the sets \( \title{F}_\alpha \) define a \(K\)-partition of \(X(\tau)\). They are clearly disjoint and their sum is whole of \(X(\tau)\).
        
        Let \( A \subseteq \kappa \) be such that \( F_A = \bigcup\limits_{\alpha \in A} F_\alpha \) is meager, ie \( F_A = \bigcup\limits_{i \in \omega} M_i \), where \(M_i\) are nowhere dense in \(X\). Let
        \[
            G_k = \{ x \in X(\tau) \colon \exists_{i \in \omega} x_i \subseteq X \setminus M_k \}.
        \]
        The sets \(G_k\) are sums of base open sets in \(X(\tau)\) and thus are open. Let \( U_0, \ldots, U_n \subseteq X \) be open and satisfying \( U = U_0 \cap \ldots \cap U_n \not= \emptyset \). As \(M_k\) is nowhere dense we have \( U \setminus \Cl(M_k) \not= \emptyset \). Then for \( G = \{ x \in X(\tau) \colon x_0 = U_0, \ldots, x_n = U_n \} \) we have \( x = (U_0, \ldots, U_n, U, U, U, \ldots) \in G \cap G_k \) and it follows that \( G_k \) is dense. Thus \( \tilde{F}_A \subseteq \{ x \in X(\tau) \colon \bigcap x \cap F_A \not= \emptyset \} \subseteq \bigcap\limits_{i \omega} (X(\tau) \setminus G_i) \) is meager.
        
        Let in turn \( A \subseteq \kappa \) be such that \( F_A = \bigcup\limits_{\alpha \in A} F_\alpha \) is non-meager. We know that \( F_A \) has the Baire property, ie \( F_A = U \triangle \bigcup\limits_{i \in \omega} M_i \), where \(M_i\) are nowhere dense in \(X\). Moreover they can be assumed to be closed. Let once more
        \[
            G_k = \{ x \in X(\tau) \colon \exists_{i \in \omega} x_i \subseteq X \setminus M_k \}
        \]
        and
        \[
            G_U = \{ x \in X(\tau) \colon \exists_{i \in \omega} x_i \subseteq U  \}.
        \]
        We already know that the sets \(X(\tau) \setminus G_k\) are nowhere dense. Let \( x \in G_U \setminus \bigcup\limits_{i \in \omega} (X(\tau) \setminus G_k) \). For some \( n \in \omega \) we have \( x_n \subseteq U \), so \( \bigcap x \subseteq U \). On the other hand for all \( n, i \in \omega \) we have \( x_n \cap M_k = \emptyset \), ie \( \bigcap x \cap M_k = \emptyset \). It follows that
        \[
            \bigcap x \subseteq U \setminus \bigcup\limits_{i \in \omega} M_i \subseteq U \triangle \bigcup\limits_{i \in \omega} M_i \subseteq F_A
        \]
        which in turn shows that \( x \in \tilde{F}_A \) and thus \( G_U \setminus \bigcup\limits_{i \in \omega} (X(\tau) \setminus G_k) \subseteq \tilde{F}_A \).
    \end{Prf}

    This result can be further refined in a case when we have a \(K\)-partition of a compact space.

    \begin{Thm}
        Let \(X\) be a compact space with \( K \)-partition \( \mathcal{F} \) of minimal cardinality \( \kappa \) and let \(\tau\) be the topology of \(X\). Then the space \(X^*(\tau)\) is a complete metric space and has a \(K\)-partition.
    \end{Thm}
    \begin{Prf}
        First we will show that \(X^*(\tau)\) is complete. Let \( (\xi_n)_{n \in \omega} \) be a Cauchy sequence in \(X^*(\tau)\). It is also a Cauchy sequence in a complete space \( \tau^\omega \) and thus there exist \( x_n \in X \) and \( N_n \in \omega \) such that for \( k \ge N_n \) we have \( \xi_{n,k} = x_n \). We need to show that \( x = (x_n)_{n \in \omega} \in X^*(\tau) \) It is clear that we can assume that \( N_{n+1} > N_n \). With that assumption we have
        \[
            \Cl(x_{n+1}) = \Cl(\xi_{n+1, N_{n+1}}) \subseteq \xi_{n, N_{n+1}} = \xi_{n, N_n} = x_n.
        \]
        This proves our claim.
        
        We define a \(K\)-partition on \( X^*(\tau) \) in the exact same way as we did for \(X(\tau)\), namely:
        \[
            \tilde{F}_\alpha = \{ x \in X(\tau) \colon \alpha = \min \{ \beta < \kappa \colon \bigcap\limits_{k < \omega} x(k) \cap F_\beta \not= \emptyset \} \} \text{ for } \alpha \in \kappa.
        \]
        For the remainder of the proof, the same reasoning as for \(X(\tau)\) is valid for \(X^*(\tau)\).
    \end{Prf}

    \section{K-partitions, precipitous ideals and measurable cardinals}
    
    As it turns out the existence of \(K\)-partitions directly implies the existence of everywhere precipitous ideals.
    \begin{Thm}
        Let \( X \) be a space with \( K \)-partition \( \mathcal{F} \) of minimal cardinality \( \kappa \). Then there exists an open non-meager subset \( U \) of \( X \) such that \( I_{\mathcal{F}|_U} \) is everywhere precipitous.
    \end{Thm}
    \begin{Prf}
        By Banach Localization Theorem there exists an open subset \( U \) such that it has no non-empty meager open subsets. Suppose \( O \) is open and non-meager in \( U \) and \( I_{\mathcal{F}|_O} \) is not precipitous. Then there exists \( S \in I^+_{\mathcal{F}|_O} \) and descending chain of functionals
        \begin{equation}
        \Phi_0 > \Phi_1 > \ldots
        \end{equation}
        on \( S \). We have \( \bigcup\limits_{\alpha \in S} F_\alpha \cap O = V \triangle M \), where \( V \) is open and \( M \) is meager. Let \( W_i = W_{\Phi_i} = \{ \dom(\varphi) \colon \varphi \in \Phi_i \} \) be corresponding \( I_{\mathcal{F}|_O} \)-partitions. Note that they are also \( I_{\mathcal{F}|_V} \)-partitions.
        
        Let \( Y \in W_i \). Then there exists \( \varphi_Y \in \Phi_i \) such that \( \dom(\varphi_Y) = Y \). We have \( \bigcup\limits_{\alpha \in Y} F_\alpha \cap V = V_Y \triangle M_Y \), where \( V_Y \) is open and \( M_Y \) is meager. Note that the sets \( V_Y \) are pairwise disjoint. Indeed if it was not the case their intersection would be open and non-empty and thus non-meager by our assumption. Let \( f_Y \colon V_Y \to \kappa \) be given by
        \begin{equation}
        f_Y(x) = \alpha \mbox{ for } x \in F_\alpha.
        \end{equation}
        By maximality of \( W_i \) the sets \( V_i = \bigcup\limits_{Y \in W_i} V_Y \) are open and dense in \( V \) and therefore by Baire theorem \( \bigcap\limits_{i \in \omega} V_i \not= \emptyset \).
        
        Let \( f_i \colon V_i \to \kappa \) be given by \( f_i = \bigcup\limits_{Y \in W_i} f_Y \). Then by the properties of functionals we have that
        \[
        \forall_{i \in \omega} \forall_{x \in V_{i+1}} f_i(x) > f_{i+1}(x),
        \]
        that is
        \[
        \forall_{i \in \omega} \forall_{x \in V_{i+1}} f_i(x) \ni f_{i+1}(x).
        \]
        Take \( x \in \bigcap\limits_{i \in \omega} V_i \). Then we have
        \[
        f_0(x) \ni f_1(x) \ni \ldots
        \]
        which contradicts the axioms of regularity.
    \end{Prf}


    
    By \cite{KunFr} we already know that the existence of \(K\)-partitions and existence of measurable cardinals are equiconsistent. What's more with some minor additional assumptions the existence of a \(K\)-partitions implies the existence of a measurable cardinal.

    \begin{Prp}
        Let \( X \) be a space with \( K \)-partition \( \mathcal{F} \). If \( \sat(I_\mathcal{F}) < \omega \) then there exists an open non-meager set \( U \) of \( X \) such that \( I_{\mathcal{F}|_U} \) is maximal.
    \end{Prp}
    
    \begin{Prf}
        If \( \sat(I_\mathcal{F}) < \omega \) then the quotient algebra \( \kappa / I_\mathcal{F} \) is finite and thus it has atoms \( [A_1], \ldots, [A_n] \). Let \( B_i = \bigcup\limits_{\alpha \in A_i} F_\alpha \). By \( K \)-partition property \( B_1 = U \triangle M \), where \( U \) is open and \( M \) is meager. Consequently \( U \cap B_i \) is meager for \( i > 1 \) and therefore \( I_{\mathcal{F}|_U} \) is maximal. 
    \end{Prf}

    \begin{Thm}
        Let \( X \) be a space with \( K \)-partition \( \mathcal{F} \) of minimal cardinality \( \kappa \). Let \( X(I_{\mathcal{F}}) \) be complete. Then \( \kappa \) is a measurable cardinal.
    \end{Thm}

    \begin{Prf}
        From minimality with respect to the condition above \( \kappa \) is regular. We will show that there exists an open subset \( U \) of \( X \) such that \( I_{\mathcal{F}|_U} \) is a maximal ideal. 
        
        Suppose that \( \sat(I_\mathcal{F}) \ge \omega_1 \). There exists
        \[ \mathcal{A} = \{ A_n \in I_{\mathcal{F}}^+ \colon n \in \omega \} \]
        such that \(A_n \cap A_m \in I_{\mathcal{F}} \) for \( n \not= m \). As \(I_{\mathcal{F}}\) is \(\omega_1\)-additive we can replace \(A_n\) with \(A_n \setminus \bigcup\limits_{m \not= n} A_m \) to get \( A_n \cap A_m = \emptyset \) for \( n \not= m \). Let \( B_n = \bigcup\limits_{k \ge n} A_k \).
        
        Now let \( x_m(n) = B_n \) for \(n \le m\) and \( x_m(n) = B_m \) otherwise. Of course all \( x_m \in X(I_\mathcal{F}) \) and \( (x_m)_{m \in \omega} \) is a Cauchy sequence, but its limit is not in \( X(I_\mathcal{F}) \) contradicting its completeness.
        
        Thus \( \sat(I_\mathcal{F}) \) is finite and by the second of the above propositions there exists and open set \( U \) in \( X \) such that \( I_{\mathcal{F}|_U} \) is maximal. By the theorem above we may assume that it is also everywhere precipitous. By precipitousness it is also \( \kappa \)-complete and thus it makes \( \kappa \) into a measurable cardinal.
    \end{Prf}

    \section{Possible connections with real-measurable cardinals and further developments}

    The are reasons to believe that \(K\)-partitions may be connected with real-measurable cardinals as well. If there exists a real-measurable cardinal \( \kappa \le \con \) then we know from \cite{Ul} and \cite{Sol} that there exists a \( \kappa \)-additive measure \(\mu\) defined on \( \Pow(\R) \) that extends the regular Lebesgue measure. Let us restrict that measure to the compact interval \([0;1]\) and define a Boolean algebra
    \[
        B(\mu) = \Pow([0;1]) / \Delta_\mu
    \]
    where \( \Delta = \{ A \subseteq [0;1] \colon \mu(A) = 0 \} \). We might now consider a Stone space \( X(\mu) = St(B(\mu)) \) and introduce a partition of \(X(\mu)\) in a following way. Let \( F_x \) be the filter on \(B(\mu)\) generated by the elements \([(x - \frac{1}{n}; x + \frac{1}{n}) \cap [0;1]]\) and \( \tilde{F}_x \) be the family of all ultrafilters extending \(F_x\). It's fairly easy to see that these families are disjoint, nowhere dense and cover all of \(X(\mu)\). Showing that they in fact do form a \(K\)-partition is planned to be a subject of future works. Along with that connections between the spaces \(X(\mu)\) and density topologies will be shown.

\end{document}